\documentclass[runningheads,a4paper]{llncs}

\usepackage{amssymb}
\usepackage{amsmath}
\usepackage{graphicx}
\usepackage{cite}

\usepackage{url}
\urldef{\mailsa}\path|{garciabmarta@uniovi.es, colubi@uniovi.es}|
\urldef{\mailsb}\path|erricos@cut.ac.cy|

\begin{document}


\title{Lasso estimation of an interval-valued multiple regression model}

\titlerunning{Lasso estimation of an  interval-valued multiple regression model}

%
%
\author{Marta Garc\'{i}a B\'{a}rzana\inst{1}
\and Ana Colubi\inst{1} \and Erricos J. Kontoghiorghes\inst{2}}
\authorrunning{Lasso estimation of an  interval-valued multiple regression model}

\institute{Department of Statistics. University of Oviedo\\
C/ Calvo Sotelo s/n, Oviedo, 33007, Spain \\
\and Department of Commerce, Finance and Shipping \\
Cyprus University of Technology. P.O. Box 50329, CY-3603 Limassol, Cyprus
\mailsa\\
\mailsb\\}

%
%

\toctitle{Lasso estimation of an  interval-valued multiple regression model}
\maketitle

\begin{abstract}
A multiple interval-valued linear regression model considering all the cross-relationships between the mids and spreads of the intervals has been introduced recently. A least-squares estimation of the regression parameters has been carried out by transforming a quadratic optimization problem with inequality constraints into a linear complementary problem and using Lemke's algorithm to solve it. Due to the irrelevance of certain cross-relationships, an alternative estimation process, the Lasso, is developed. A comparative study showing the differences between the proposed estimators is provided.
\end{abstract}

\section{Introduction}\label{MartaGB_sec:1}
Intervals represent a powerful tool to capture the imprecission of certain characteristics that cannot be fully described with a real number. For example, the measures provided by instruments which have some errors in their measurements \cite{Abdallah:07}. Moreover, intervals also model some features which are inherently interval-valued. For instance, the range of variation of the blood preasure of a patient along a day \cite{Blancoetal:11} or the tidal fluctuation \cite{Ramos:13}.
 
The statistical study of regression models for interval data has been extensively addressed lately in the literature  \cite{Blancoetal:11,Blancoetal:14,Giordani:11,GRetal:07,LimaNetoetal:10}, deriving into several alternatives to tackle this problem. On one hand, the estimators proposed in \cite{Giordani:11,LimaNetoetal:10} account the non-negativity constraints satisfied by the spread variables, but do not assure the existence of the residuals. Hence, they can lead to ill-defined estimated models. On the other hand, the models proposed in \cite{Blancoetal:11,Blancoetal:14,GRetal:07} are formalized according to the natural interval arithmetic and their estimators lead to models that are always well-defined over the sample range.

The multiple linear regression model \cite{Blancoetal:14} considered belongs to the latter approach and its main advantage is the flexibility derived from its way to split the regressors, allowing us to account for all the cross-relationships between the centers and the radious of the interval-valued variables. Nevertheless, this fact entails an increase in the number of regression parameters and thus, a Lasso estimation is considered in order to schrink some of these coefficients towards zero. The Lasso estimation of an interval-valued regression model has been previously addressed in \cite{Giordani:11}, but being this a more restrictive model formalized in the first framework.

The paper is organized as follows. Section
\ref{MartaGB_Preliminaries:2} presents some preliminary concepts about the
interval framework and Section \ref{MartaGB_model:3} contains the formalization of the model. The Least-Squares and Lasso estimations of the proposed model are developed in subsections \ref{MartaGB_LS:3} and \ref{MartaGB_Lasso:3}. Section \ref{MartaGB_Gio:4} briefly describes the Lasso model proposed by Giordani \cite{Giordani:11}. The empirical performance of the estimators proposed in Sections \ref{MartaGB_model:3} and \ref{MartaGB_Gio:4} is compared in Section \ref{MartaGB_Example} by means of a illustrative real-life example.
Section \ref{MartaGB_Conclusions:6} finishes with some conclusions.

\section{Preliminaries}\label{MartaGB_Preliminaries:2}
\noindent Interval data are defined as elements belonging to the space $\mathcal{K_{\mathrm{c}}}(\mathbb{R})=\{[a_1,a_2] : a_1,a_2 \in \mathbb{R}, a_1
\le a_2\}$. Given an interval $A\in \mathcal{K_{\mathrm{c}}}(\mathbb{R})$, it can be parametrized in terms of its center or {\it midpoint}, $\textrm{mid}\ \hspace*{-0.1cm}  A =(\sup A + \inf A)/2$, and its radius or {\it spread}, $\textrm{spr}\ \hspace*{-0.1cm}  A = (\sup A - \inf A)/2$. Nonetheless, intervals can alternatively be expressed by means of the so-called canonical decomposition \cite{Blancoetal:11} given
by $A=\mathrm{mid}A[1\pm 0] + \mathrm{spr}A[0\pm 1]$. This decomposition allows us to consider separately the {\it mid} and {\it spr} components of $A$, which will lead into a more flexible model. The interval arithmetic on
$\mathcal{K _{\mathrm{c}}}(\mathbb{R})$ consists of the Minkowski addition and the product by scalars defined as follows by the jointly expression: 
$A+ \delta B=[( \mathrm{mid}A + \delta\mathrm{mid}B ) \ \pm \
(\mathrm{spr} A + |\delta| \ \mathrm{spr} B)]$ for any $A,B \in
\mathcal{K _{\mathrm{c}}}(\mathbb{R})$ and $\delta \in \mathbb{R}$.

The space $(\mathcal{K _{\mathrm{c}}}(\mathbb{R}),+,\cdotp)$ is not linear but semilinear, as the existence of symmetric element with respect to the addition is not guaranteed in general. An additional operation is introduced, the so-called Hukuhara difference between the intervals $A$ and $B$. The difference $C$ is defined as $C=A-_H B \in \mathcal{K _{\mathrm{c}}}(\mathbb{R})$ verifying that $A=B+C$. The existence of $C$ is subject to the fulfillement of the expression $\mathrm{spr} B \le \mathrm{spr} A$.

\noindent Given the intervals $A,B\in \mathcal{K _{\mathrm{c}}}(\mathbb{R})$, the metric $d_{\tau}(A,B)=((1-\tau)\,((\mathrm{mid}A-\mathrm{mid}B)^2+\tau\,(\mathrm{spr}A-\mathrm{
spr}B)^2))^{\frac{1}{2}}$, for an arbitrary $\tau \in (0,1)$, is the $L_2$-type distance to be considered. $d_{\tau}$ is based on the metric $d_{\theta}$ defined in \cite{Trutschnigetal:09}. 

\noindent Given a probability space $(\Omega, \mathcal {A}, P)$ the mapping $\emph{\textbf{x}}:\Omega \rightarrow \mathcal{K _{\mathrm{c}}}(\mathbb{R})$ is a random interval iff it is a measurable Borel mapping. The moments to be considered are the classical Aumann expected value for intervals; the variance defined following the usual Fr\'{e}chet variance \cite{Nather:97} associated with the Aumann expectation in the interval space $(\mathcal{K_{\mathrm{c}}}(\mathbb{R}), d_\tau)$; and the covariance defined in terms of mids and spreads as $\sigma_{\emph{\textbf{x}},\emph{\textbf{y}}}=(1-\tau)\,\sigma_{\textrm{mid}
\emph{\textbf{x}},\textrm{mid} \emph{\textbf{y}} }+ \tau\sigma_{\textrm{spr}
\emph{\textbf{x}},\textrm{spr} \emph{\textbf{y}}}$.

\section{The multiple linear regression model}\label{MartaGB_model:3}
 Let \emph{\textbf{y}} be a response random interval and let $\emph{\textbf{x}}_1,\emph{\textbf{x}}_2,\ldots,\emph{\textbf{x}}_k$ be $k$ explanatory random intervals. The model is formalized in a matrix notation as follows:
\begin{equation}\label{modelmatrix}
\emph{\textbf{y}}=X^{Bl}\,B+\boldsymbol{\varepsilon} \ ,
\end{equation}
where $B=(b_1|b_2|b_3|b_4)^t\in \mathbb{R}^{4k \times 1}$ with $b_i\in \mathbb{R}^k$ ($i\in \{1,2,3,4\}$), $X^{Bl}=(\boldsymbol{x^M|x^S|x^C|x^R})\in\mathcal{K_{\mathrm{c}}}(\mathbb{R})^{1\times
4k}$ where the elements are defined as $\boldsymbol{x^M}=\mathrm{mid}\,x^t\,[1\pm 0]$, $\boldsymbol{x^S}=\mathrm{spr}\,x^t\,[0\pm 1]$, $\boldsymbol{x^C}=\mathrm{mid}\,x^t\,[0\pm 1]$ and $\boldsymbol{x^R}=\mathrm{spr}\,x^t\,[1\pm 0]$, considering the canonical decomposition of the regressors. \newline $\mathrm{mid}\,x=(\mathrm{mid}\,\boldsymbol{x_1},\mathrm{mid}\,\boldsymbol{x_2},\ldots,\mathrm{mid}\,\boldsymbol{x_k})^t\in\mathbb{R}^k$ (analogously $\mathrm{spr}\,x$) and $\boldsymbol{\varepsilon}$ is a  random interval-valued error such that $E(\boldsymbol{\varepsilon}|x)=\Delta \in
\mathcal{K_{\mathrm{c}}}(\mathbb{R})$.

The following separate linear relationships for the {\it mid} and {\it spr} components of the intervals are derived from (\ref{modelmatrix}): 
\vspace*{-0.2cm}
\begin{subequations}
\begin{equation}\label{separate-models-a}
\mathrm{mid}\,\emph{\textbf{y}}=\mathrm{mid}\,x^t\,b_1+\mathrm{spr}\,x^t\,
b_4+\mathrm{mid}\,\boldsymbol{\varepsilon}\ ,
\end{equation}
\vspace*{-1cm}

\begin{equation}\label{separate-models_b}
\mathrm{spr}\,\emph{\textbf{y}}=\mathrm{spr}\,x^t\,|b_2|+|\mathrm{mid}\,x^t|\,
|b_3|+\mathrm{spr}\,\boldsymbol{\varepsilon}\ .
\end{equation}
\end{subequations}

\noindent Thus, the flexibility of the model arises from the possibility of considering all the information provided by $\mathrm{mid} \emph{x}$ and $\mathrm{spr} \emph{x}$ to model $\mathrm{mid} \emph{\textbf{y}}$ and
$\mathrm{spr} \emph{\textbf{y}}$, as follows from (\ref{separate-models-a}) and (\ref{separate-models_b}). This represents an improvement with respect to previous models that merely addressed the relationship between the mids of the variables or between the spreads but never any cross-relationship (mid-spr).

Nevertheless, the inclusion of more coefficients entails an increase in the dimensionality of the estimation process. Some of these coefficients could be zero as not all the new introduced variables will contribute. Therefore it is proposed to estimate (\ref{modelmatrix}) by least-squares and by Lasso and compare the advantages and disadvantages that each estimation process provides.

\subsection{The Least-Squares estimation}\label{MartaGB_LS:3}
Given $\{\left(\boldsymbol{y_j},\boldsymbol{x_{i,j}}\right) :
i=1,\ldots,k, j =1,\ldots,n\}$ a simple random sample of intervals obtained from $(\boldsymbol{y},\boldsymbol{x_1},\ldots,\boldsymbol{x_k})$ in
(\ref{modelmatrix}) the estimated model is
\begin{equation}\label{samplemodel}
\widehat{y}=X^{ebl}\widehat{B}+\widehat{\varepsilon} 
\end{equation}
\noindent where $y=(\emph{\textbf{y}}_1,\ldots,\emph{\textbf{y}}_n)^t$, $X^{ebl}=(X^M|X^S|X^C|X^R)\in \mathcal{K_{\mathrm{c}}}(\mathbb{R})^{n\times 4k}$, $\varepsilon=(\boldsymbol\varepsilon_{1},\ldots,\boldsymbol\varepsilon_{n})^t$ is such that $E(\varepsilon|x)=1^n \Delta$ and $B$ as in (\ref{modelmatrix}).
$X^M$ is the $(n\times k)$-interval-valued matrix such that $(X^M)_{j,i}=\mathrm{mid}\emph{\textbf{x}}_{i,j} [1\pm 0]$ (analogously $X^S,X^C$ and $X^R$). Given an arbitrary vector of regression coefficients $A\in \mathbb{R}^{4k\times 1}$ and an interval of residuals $C\in \mathcal{K}_c(\mathbb{R})$, the Least-Squares estimation looks for $\widehat{B}$ and $\widehat{\Delta}$ minimizing the distance $\displaystyle d_{\tau}^2(y, X^{ebl}A+1^n C)$. $\widehat{\Delta}$ can be obtained separately and firstly by the expression $ \widehat{\Delta}= \overline{\emph{\textbf{y}}}-_H
\overline{X^{ebl}}\widehat{B}.$

Recalling that, by definition, $X^S=-X^S$ (and analogously $X^C=-X^C$) the estimation process of the coefficients $b_2$ and $b_3$ accompanying these variables can be simplified by searching only for non-negative estimates. By contrast, coefficients $b_1$ and $b_4$ are not affected by any kind of restrictions so they can be estimated directly by OLS. Moreover, it has to be assured the existence of the residuals defined as the Hukuhara differences $\varepsilon= \emph{\textbf{y}} -_H X^{ebl}B$. For this purpose the minimization problem ends up to be the following constrained quadratic problem:
\vspace*{-0.3cm}
\begin{eqnarray}\label{quadratic-objective-function}
\min_{A_m\in\ \mathbb{R}^{2k}, \ A_s\in \ \Gamma}(1-\tau)\,\Vert v_m-F_m A_m\Vert^2+\tau\,\Vert v_s-F_s A_s\Vert^2 \\
\nonumber \Gamma=\{(a_2,a_3)\in
[0,\infty)^k\times[0,\infty)^k :
\mathrm{spr}\,X\,a_2+|\mathrm{mid}\,X|\,a_3\leq \mathrm{spr}\,y\} ,
\end{eqnarray}
\vspace*{-0.5cm}

\noindent being
$v_m=\mathrm{mid}y-\overline{\mathrm{mid}\boldsymbol{y}}\,1^n$,
$v_s=\mathrm{spr}y-\overline{\mathrm{spr}\boldsymbol{y}}\,1^n\in
\mathbb{R}^n$,
$F_m=\mathrm{mid}X^{ebl}-1^n(\overline{\mathrm{mid}X^{ebl}})$,
$F_s=\mathrm{spr}X^{ebl}-1^n(\overline{\mathrm{spr}X^{ebl}})\in\mathbb{R}^{n\times
2k}$, $A_m=(a_1 | a_4)^t\in \mathbb{R}^{2k\times 1}$ the
coefficients related to the midpoints and $A_s=(a_2 |
a_3)^t\in\mathbb{R}^{2k\times 1}$ the coefficients related to the
spreads, with $a_l\in \mathbb{R}^k$, $l=1,\ldots,4$.

There are several numerical ways to tackle the resolution of a quadratic problem as (\ref{quadratic-objective-function}). Given the shape of the objective function, the minimization process is solved separately over $A_m$ and $A_s$. Those coefficients related with the mids ($A_m$) are not affected by constraints and therefore, the OLS estimator can be used directly. Thus $\widehat{A_m}=(F_m^tF_m)^{-1}F_m^t v_m$. However, in order to proceed with the constrained minimization over $A_s$, Karush-Kuhn-Tucker conditions guarantee the existence of local optima solution, which can be computed with standard numerical tool. Nevertheless, in order to obtain an exact solution and a more handy estimator of $A_s$, (\ref{quadratic-objective-function}) can be equivalent expressed as a {\it Linear Complementary Problem} with the shape:
\begin{eqnarray}\label{LCP}
  \omega=M\,\lambda+q \hspace*{0.5cm}  \textit{s.t.\ } \hspace*{0.3cm} \omega,\lambda\ge 0 \  , \  \omega_j\lambda_j=0 \  , \ j=1,\ldots, n+1\ ,
\end{eqnarray}
\noindent with $M=(R\,Q^{-1}\,R^t)$ and $q=(-R\,Q^{-1}\,c-r)$ (details in \cite{Blancoetal:14}).
Thereby, once $\lambda$ is obtained, the expression of the estimator is $\widehat{A_s}=Q^{-1}\,(R^t\,\lambda-c)$.

\subsection{The Lasso estimation}\label{MartaGB_Lasso:3}
Lasso, {\it Least Absolute Shrinkage and Selection Operator}, is a regression method that involves penalizing the sum of the absolute values of the regression coefficients estimates. For this purpose it involves a regularization parameter which affects directly the estimates: the larger the value of this parameter, the more estimates that are shrunk towards zero. This coefficient cannot be estimated statistically, so a cross-validation process is usually applied. 

As previously, (\ref{quadratic-objective-function}) can solved separately. On one hand, the classical Lasso method will be used to obtain the estimator of the regression coefficients related to the mids. Then, the problem is expressed as:
\vspace*{-0.3cm}
\begin{eqnarray}
\nonumber \frac{1}{2}\Vert v_m-A_m\,F_m\Vert_2^2+\lambda\, \sum_{j=1}^{2k}|A_{m_j}| 
\end{eqnarray}
\vspace*{-0.6cm}

 \noindent being $\lambda$ the regularization parameter. There are different programs capable to solve this problem (such as Matlab or R). The lasso.m Matlab function is the one used to obtain $\widehat{A_m}$.

On the other hand, for those coefficients related with the spreads a constrained Lasso algorithm has been developed as a modified version of the code proposed by Mark Schmidt (2005) \cite{Schmidt} and is available upon request. The problem is given by:
\vspace*{-0.3cm}
\begin{eqnarray}
\nonumber \frac{1}{2}\Vert v_s-A_s\,F_s\Vert_2^2 +\lambda\sum_{j=1}^{2k}|A_{s_j}| \hspace*{0.5cm}  {\rm s.t} \hspace*{0.2cm} R A_s\ge r.& & 
\end{eqnarray}
\vspace*{-0.6cm}

The most usual elections of $\lambda$ are the value than minimizes the Cross-Validation Mean Square Error ($\lambda_{MSE}$) and the value that provides a simpler or more parsimonious model with respect to $\lambda_{MSE}$ (in terms of more zero coefficients) but at the same time with one-standard-error ($\lambda_{1SE}$).

\section{Giordani's Lasso estimation}\label{MartaGB_Gio:4}
The so-called {\it Lasso-based Interval-valued Regression (Lasso-IR)} proposed by Giordani in \cite{Giordani:11} is another Lasso method to deal with a multiple linear regression model for interval data. However, the later regression model is not formalized following the interval arithmetic and can end up with an ill-defined estimated model. Keeping the same notation as in (\ref{separate-models_b}), it requires the non-negativity of $b_2$ and $b_3$ but does not test if the Hukuhara's difference $\varepsilon=\emph{\textbf{y}} -_H X^{ebl}B$ exists. The optimization problem can be written (analogously to (\ref{quadratic-objective-function})) as:
\vspace*{-0.3cm}
\begin{eqnarray}\label{MartaGB-modelo-Gio}
\min_{A_m,As} (1-\tau)\,\Vert v_m-F_m\,A_m\Vert^2+\tau \Vert v_s-F_s(A_m+A_a)\Vert^2\\
\nonumber \hspace*{-2cm} F_s (A_m+A_a)\ge 0, \sum_{j=0}^p|A_{a_j}|\le t
\end{eqnarray}
\vspace*{-0.4cm}

The coefficients related to the spreads ($A_s$) are the ones for the mids ($A_m$) plus a vector of additive coefficients ($A_a$) showing the distance that they are allowed to differ from $A_m$. In this case (\ref{MartaGB-modelo-Gio}) has been expressed as a constrained quadratic problem, where there is a one-to-one correspondence between $\lambda$ and $t$. The value of $t$ that minimizes the cross-validation mean square error is the one considered. In order to solve the problem a stepwise algorithm based on \cite{Lawson:95} is proposed.

 Another important difference, which entails less flexibility in the model, is the limitation of being able to study separately the relationships between the mids and the relationship between the spreads of the intervals but never any cross-relationship.

\begin{remark}\label{Remark1}
There is a particular case of model (\ref{MartaGB_sec:1}), the so-called Model $M$ addressed in \cite{Blancoetal:11}, which is formalized in the interval framework but has the same lack of flexibility as (\ref{MartaGB-modelo-Gio}). In this case $b_3$ and $b_4=(0,\ldots,0)$, so the model has the shape:
\begin{equation}\label{MartaGB-Modelo-M}
\boldsymbol{y}=b_1\,\boldsymbol{x^M}+b_2\,\boldsymbol{x^S}+\boldsymbol{\varepsilon}.
\end {equation} 
\end{remark}

\section{A real-life illustrative example}\label{MartaGB_Example}
The following example contains the information of a sample of 59 patients (from a population of 3000) hospitalized in the Hospital Valle del Nal\'{o}n in Asturias, Spain. The variables to be considered are the ranges of fluctuation of the diastolic blood preasure over the day ($\boldsymbol{y}$), the pulse rate ($\boldsymbol{x_1}$) and the systolic blood preasure ($\boldsymbol{x_2}$). The dataset can be found in \cite{Blancoetal:11} and \cite{GRetal:07}.

In order to make possible the comparison between the estimator proposed in Sect. \ref{MartaGB_Gio:4} and those ones introduced in Subsect. \ref{MartaGB_LS:3} and Subsect. \ref{MartaGB_Lasso:3}, the example will be developed for the simpler model explained in Remark \ref{Remark1}. 

Given the displayed model in (\ref{MartaGB-Modelo-M}), $\boldsymbol{y}=b_{1}\boldsymbol{x^M_1}+b_{2}\boldsymbol{x^M_2}+b_{3}\boldsymbol{x^S_1}+b_{4}\boldsymbol{x^S_2}+\varepsilon$, the estimates of the regression coefficients are summarized in Table 1:
\begin{table}
\caption{Estimates of the regression coefficients for the three estimators: LS, Lasso (for the two more representatives values of $\lambda$) and Lasso-IR (for a fixed value of $t$=0.10 prefixed by the author). The last column contains the MSE of the models mimicking its definition in the classical framework.}
\label{MartaGB_table:1}       
%
%
\begin{tabular}{l c c c c c}
\hline\noalign{\smallskip}
  &  $\widehat{b_{1}}$ &  $\widehat{b_{2}}$ &  $\widehat{b_{3}}$& $\widehat{b_{4}}$ &MSE  \\
\hline
$LS-estimation \,(Sect. \,\ref{MartaGB_LS:3})$ &0.4497 &0.0517& 0.2588&0.1685& 68.2072\\
\hline
$Lasso-estimation \,(Sect. \,\ref{MartaGB_Lasso:3})$& 0.4202&0.0020& 0.3379&0.2189 & 68.8477\\
 $\tiny{ \lambda_{MSE} }$& \tiny{ (0.6094) }& &\tiny{ ( 0.0259) }& &\\
\hline
$Lasso-estimation \,(Sect. \,\ref{MartaGB_Lasso:3})$& 0.2749& 0& 0.0815 &0&  76.9950\\
$\tiny{ \lambda_{1SE} }$&  \,\tiny{ (3.2521) }& &\tiny{ (1.8736) }& &\\
\hline
$Lasso-IR \,(Sect. \ref{MartaGB_Gio:4})$&  0.5038&0.1261& 0.4847 & 0.3605& 71.2418\\
\hline
\end{tabular}
\end{table}

In view of the results in Table \ref{MartaGB_table:1}, those coefficients which take small values with the LS-estimation ($\widehat{b_2}$ and $\widehat{b_4}$) are schrunk towards zero with the most preferable Lasso estimation (for $\lambda_{1SE}$). However, this entails a significant increase of the MSE. In the case of using our Lasso-estimator for $\lambda_{MSE}$, the MSE is smaller but it does not provide a parsimonious model, being therefore its usefulness questionable. The estimator proposed in \cite{Giordani:11} reaches a high value of MSE (worse in comparison with the lasso for $\lambda_{MSE}$) and does not end up with an easy-to-interpret model.

\section{Conclusions}\label{MartaGB_Conclusions:6}
On one hand, a recently studied regression model for interval data, allowing to study all the cross-relationships between the mids and spreads of the interval-valued variables involved, is considered. This flexibility derives into an increase of the dimensionality of the model. Therefore a Lasso estimation seems appropiate to tackle this problem by setting some of these coefficients to zero. Nonetheless, a comparisson study gathering the double estimation process conducted (first by Least-Squares and after by Lasso) is provided. 

On the other hand, it is considered the Lasso-based interval-valued regression model (Lasso-IR) proposed in \cite{Giordani:11}. This model is not constrained to guarantee the existance of the residuals so it can provide misleading estimations. Moreover, it has a lack of flexibility as it solely tackles the relastionships of type mid-mid and spr-spr but no cross-relationship mid-spr.

A real-life example illustrating the difference between the estimators in terms of MSE and simplicity has been conducted.

\subsubsection*{Acknowledgments.} Financial support from the Spanish Ministerio de Ciencia e Innovacion (MICINN) through Ayuda Puente SV-PA-13-ECOEMP-66 and Acci\'{o}n Integrada PRI-AIBDE-2011-1197, is kindly acknowledged.

\bibliographystyle{splncs03}

\end{document}